\input amstex
\magnification=1200
\documentstyle{amsppt}
\NoRunningHeads
\NoBlackBoxes
\topmatter
\title Differential interactive games:\linebreak 
The short-term predictions\endtitle
\author Denis V. Juriev\endauthor
\affil ul.Miklukho-Maklaya 20-180, Moscow 117437 Russia\linebreak
E-mail: denis\@juriev.msk.ru\endaffil
\date math.OC/9901074\enddate
\endtopmatter
\document
The mathematical formalism of differential interactive games, which extends 
one of ordinary differential games [1] and is based on the concept of an
interactive control, was proposed by the author [2] to take into account 
the complex composition of controls of a real human person, which are often 
complicated couplings of his/her cognitive and known controls with the 
unknown subconscious behavioral reactions.

In the article [3] some special class of differential interactive games, the 
laced interactive games, was considered. Besides other results a mechanism of 
short-term predictions for processes in such games was proposed. It is based 
on some approximations of the laced interactive games by ordinary differential 
games.

The goal of this research note is to propose similar mechanism of heuristic 
short-term predictions for general differential interactive games.

\subhead 1. The differential interactive games\endsubhead

\definition{Definition 1 [2]} {\it An interactive system\/} (with $n$
{\it interactive controls\/}) is a control system with $n$ independent controls 
coupled with unknown or incompletely known feedbacks (the feedbacks, which are
called the {\it behavioral reactions}, as well as their couplings with 
controls are of a so complicated nature that their can not be described 
completely). {\it An interactive game\/} is a game with interactive controls 
of each player.
\enddefinition

Below we shall consider only deterministic and differential interactive
systems. For symplicity we suppose that $n=2$. In this case the general
interactive system may be written in the form:
$$\dot\varphi=\Phi(\varphi,u_1,u_2),\tag1$$
where $\varphi$ characterizes the state of the system and $u_i$ are
the interactive controls:
$$u_i(t)=u_i(u_i^\circ(t),\left.[\varphi(\tau)]\right|_{\tau\le t}),$$
i.e. the independent controls $u_i^\circ(t)$ coupled with the feedbacks on
$\left.[\varphi(\tau)]\right|_{\tau\le t}$. One may suppose that the
feedbacks are integrodifferential on $t$ in general, but below we shall 
consider only differential dependence. It means that 
$$u_i(t)=u_i(u_i^\circ(t),\varphi(t),\dot\varphi(t),\ddot\varphi(t),\ldots,
\varphi^{(n)}(t)).\tag 2$$
It is reasonable to suppose that all functional dependencies (1) and (2) are
smooth.

\subhead 2. Short-term predictions. Basic procedure\endsubhead

Let $u_i$ and $u^\circ_i$ ($i=1,2$) have $n$ degrees of freedom. Let us
consider $2n+1$ arbitrary functions $p_j(\vec u,\vec u^\circ,\varphi)$ of 
$\vec u=(u_1,u_2)$, $\vec u^\circ=(u^\circ_1,u^\circ_2)$ and $\varphi$ ($j=1,
2,\ldots,2n+1$). The knowledge of the processes in the game at $=\tau<t$ allows 
to consider $2n$ magnitudes $f_i(\tau)=\sum_{j=1}^{n+1}\alpha_{ij}(\tau)
p_j(\vec u(\tau),\vec u^\circ(\tau),\varphi(\tau))$ ($i=1,2,\ldots 2n$) such 
that $\dot f_i(\tau)\equiv 0$. One may suppose that the coefficients
$\alpha_{ij}(\tau)$ are continuous and, moreover, belong to the Lipschitz
class. Their differentiability is too strong condition to be satisfied in
practice. 

For the fixed moment $t$ let us consider $\Delta t>0$ such that the
Jacobi matrix of the mapping $\vec u\mapsto(f_1,\ldots f_n)$ is nondegenerate
at the moment $\tau=t-\Delta t$ and at the point $\vec u=\vec u(\tau)$. Under
these conditions one can locally express $\vec u$ via $\vec u^\circ$ and 
$\varphi$: $$\vec u(\tau)=\vec U_\tau(\vec u^\circ(\tau),\varphi(\tau);
f_1(\tau),\ldots f_{2n}(\tau)).\tag 3$$

The obtained relations may be used for an approximation of the interactive 
game by ordinary differential games. Let us consider a fixed moment $t_0$.
For $t>t_0$ the interactive controls $u_i(t)$ will be replaced by their
approximations $u^*_i(t)$ as in the evolution equations of the game as
in the expressions for the functions $f_i$. The magnitudes $u^*_i(t)$ are
defined by the formulas
$$\vec u^*(t)=\vec U_{t-\Delta t}(\vec u^\circ(t),\varphi(t-\Delta t);
f_1(t-\Delta t),\ldots f_{2n}(t-\Delta t))\tag 4$$
for $t>t_0$ (for $t<t_0$ they coincide with the interactive controls
$u_i(t)$). Note that $f_i$ were calculated by use of the values of
$\vec u^*$ at the moment $t-\Delta t$. 
Thus, we receive an ordinary differential game with retarded (delayed)
arguments, to which the more or less standard analysis of ordinary 
differential games can be applied. The approximation (4) generalizes
the retarded control approximation of the article [3]. The values of the 
state $\varphi$ calculated at the moment $t-\Delta t$ may be changed to
its values calculated at the moment $t$ or at the intermediate moment
$t-\alpha\Delta t$, where the parameter $\alpha\in[0,1]$ is chosen to
provide the best approximation.

Note that really we constructed a series of ordinary differential games
pa\-ra\-met\-rized by $t_0$. The obtained predictions may be used as 
short-term predictions for processes in the initial interactive game. 
Certainly, as it was marked in [3] it is difficult to perceive and to 
interpret the analytically represented results in real time. Thus, it is 
rather reasonable to use some visual representation for the series of the 
approximating games. Thus, we are constructing an enlargement of the 
interactive game, in which the players interactively observe the visual 
predictions for this game in real time. Of course, such enlargement may 
strongly transform the structure of interactivity of the game (i.e. to 
change the feedbacks entered into the interactive controls of players).

\subhead 3. Short-term predictions. Further developments\endsubhead

The basic procedure exposed above essentially depends on the choice of
the functions $p_j(\vec u,\vec u^\circ,\varphi)$. Its further developments
are based on the attempts to choose them dynamically in the most effective
way.

The simplest improvement is in the consideration of several sets
$\{p_j^{(\mu)}\}$ of such functions. The index $\mu$ labels an individual
set. Fixing the moment $t_0$ and $\Delta t$ one performs the basic procedure 
for each $\mu$ starting at $t_0-\Delta t$ instead of $t_0$. The obtained
short-term predictions for $\varphi$ are compared with the real data
in the time interval $t_0-\Delta t<t<t_0$ (at this interval $\vec u^\circ$
coincides with its observed value). The best prediction determines
$\mu$, which is used for the short-term predictions for $t>t_0$.

The index $\mu$ may vary over a finite set or over a smooth manifold.
For example, let us consider the set of $2n+2$ functions $p_j(\vec u,
\vec u^\circ,\varphi)$. They generate a linear space $V$. Any hyperplane
in this space is spanned by $2n+1$ functions, which may be used in the
basic procedure. In this case $\mu$ labels a hyperplane in the $2n+2$ 
dimensional space $V$ and, therefore, belongs to the $2n+1$-dimensional
projective space $\Bbb P^{2n+1}=\Bbb P(V^*)$. Dynamics in the interactive 
game determines a curve $\mu(t)$ in $\Bbb P^{2n+1}$. The point $\mu(t_0)$ is
the index of the best prediction constructed as above for the moment $t_0$.
The curve $\mu(t)$ may be discontinuous.

The next improvement is based on the dynamical selection of the considered
sets of functions $\{p_j^{(\mu)}\}$ with finite number of $\mu$ during the 
game. One uses the procedure above to construct an individual approximation 
at the fixed moment $t_0$. Let the set labelled by $\mu_0$ gives the worst 
prediction. In the moment $t_0+\Delta t$ one adds any new set to the 
considered ensemble of sets instead of the $\mu_0$-th one, repeat the 
procedure for this moment and so on. One may specify various algorithms to 
choose the new set.

\subhead 4. Conclusions\endsubhead

Thus, several heuristic procedures of the short-term predictions for 
processes in the differential interactive games were considered. Note that
the problem of an estimation of precision of such predictions is not
correct in view of interactivity of the game. One may only say that
in any finite time interval $t_0<t<t_1$ the prediction becomes heuristically 
better with $\Delta t\downarrow 0$. At least, it may be reasonably
effective only for rather short intervals $(t_0,t_1)$\footnote
{To estimate the maximal admitted $t_1$ one may use the following
procedure: let us consider two approximations started at $t-\Delta t_1$ and
$t-\Delta t_2$, the moment $t_1$ is defined as the maximally possible one
providing that the divergence of two various predictions is not too large.
\newline}. Nevertheless, in practice the interactive
effects are essential only on the short time intervals and the short-term
analysis of the interactive game strategically reduces it to an ordinary
game. The main problem here is to extract the necessary data from such
analysis to define this new game; here, the investigation of series of
approximating differential games and the unraveling of algebraic
correlations between them (in spirit of the nonlinear geometric algebra)
is apparently crucial (cf.[3]).

\Refs
\roster
\item"[1]" Isaaks R., Differential games. Wiley, New York, 1965;\newline
Owen G., Game theory, Saunders, Philadelphia, 1968.
\item"[2]" Juriev D., Interactive games and representation theory. I,II.
E-prints: math.FA/9803020, math.RT/9808098.
\item"[3]" Juriev D., The laced interactive games and their {\sl a
posteriori\/} analysis. E-print:\linebreak math.OC/9901043.
\endroster
\endRefs
\enddocument